\documentclass[12pt]{article}
\usepackage[T1]{fontenc}
\usepackage[ansinew]{inputenc}
\usepackage{amsmath}
\usepackage{graphicx}
\title{An almost-Schur type lemma for  symmetric $(2,0)$ tensors  and applications}
\author{Xu Cheng \thanks{the author is  partially supported by CNPq and Faperj  of Brazil.}}
\usepackage{amsmath, amssymb}
\usepackage{amsthm}
\makeatletter

\@addtoreset{equation}{section} \makeatother
\newtheorem{thm}{Theorem}[section]
 \newtheorem{lem}{Lemma}[section]
\newtheorem{cor}{Corollary}[section]
 \newtheorem{rema}{Remark}[section]
\newtheorem{prop}{Proposition}[section]\newtheorem{Def}{Definition}[section]

\begin{document}
\date{}
\maketitle
 \begin{abstract}
 
In our previous paper in \cite{C}, we generalized  the almost-Schur lemma of De Lellis and Topping for closed manifolds with nonnegative Rcci curvature to any closed manifolds.
In this paper, we generalize the above results to symmetric $(2,0)$-tensors and give the applications including  $r$th mean curvatures of closed hypersurfaces in a space form and $k$ scalar curvatures for closed locally conformally flat  manifolds. 
 \end{abstract}

\baselineskip=14pt
\section {Introduction}

Recall that an $n$-dimensional Riemannian manifold $(M, g)$ is called  to be Einstein if its traceless Ricci tensor $ \text{\r{Ric}}=\text{Ric}-\frac{R}{n}g$
 is identically zero. Here $\text{Ric}$ and   $R$ denote Ricci curvature and scalar curvature respectively.  Schur's lemma states that the scalar curvature of an Einstein manifold of dimension $n\geq 3$ must be constant.  
In  \cite{dLT} , De Lellis and Topping  discussed the stability and rigidity of Schur's lemma for closed manifold and  proved the following almost Schur lemma, as they called.
\begin{thm} (\cite{dLT})\label{thm-1}  If $(M, g)$ is a closed Riemannian manifold of dimension $n$ with nonnegative Ricci curvature, $ n\geq 3$, then 
\begin{equation}\label{ine-i1-1}\int_{M}(R-\overline{R})^2\leq \frac{4n(n-1)}{(n-2)^2}\int_{M}|\text{Ric}-\frac{R}{n}g|^2,
\end{equation}
and equivalently, 
\begin{align}\label{ine-i1} \int_{M}|\text{Ric}-\frac{\overline{R}}{n}g|^2\leq \frac{n^2}{(n-2)^2}\int_{M}|\text{Ric}-\frac{R}{n}g|^2,
\end{align}
where $\overline{R}=\frac{\int_MRdv}{\text{vol}(M)}$ denotes the average of $R$ over $M$. Moreover the equality in (\ref{ine-i1-1}) or (\ref{ine-i1}) holds if and only if $M$ is Einstein.
\end{thm}
De Lellis and Topping \cite{dLT} also proved their estimates are sharp.  First,  the constants are optimal in (\ref{ine-i1-1}) and (\ref{ine-i1}) (\cite{dLT}, Section $2$). Second, the curvature condition $\text{Ric}\geq 0$ cannot simply be dropped (see the examples in the proof of Prop. 2.1 and 2.2 in \cite{dLT}).  Without the condition of nonnegativity of Ricci curvature,   the same type of inequalities  as (\ref{ine-i1-1}) and (\ref{ine-i1}) cannot hold if  the constants in these  inequalities  only depend on the lower bound of the Ricci curvature.
In \cite{C}, we considered  the case of closed manifolds without the hypothesis of nonnegativity of Ricci curvature and proved that
\begin{thm} (\cite{C}) \label{thm-2}  If $(M, g)$ is a closed Riemannian manifold of dimension $n, n\geq 3,$ then 
\begin{equation}\label{ine-i-1}\int_{M}(R-\overline{R})^2\leq \frac{4n(n-1)}{(n-2)^2}\left(1+\frac{nK}{\lambda_1}\right)\int_{M}|\text{Ric}-\frac{R}{n}g|^2,
\end{equation}
and equivalently,
\begin{equation}\label{ine-i-2} \int_{M}|\text{Ric}-\frac{\overline{R}}{n}g|^2\leq \frac{n^2}{(n-2)^2}\left[1+\frac{4(n-1)K}{n\lambda_1}\right]\int_{M}|\text{Ric}-\frac{R}{n}g|^2,
\end{equation}
where  $\lambda_1$ denotes the first nonzero eigenvalue of Laplace operator on $(M,g)$,
 $K$ is nonnegative constant such that  the Ricci curvature of $(M,g)$ satisfies $\text{Ric}\geq -(n-1)K$, and $\overline{R}$ denotes the average of $R=\frac{\int_MRdv}{\text{vol}(M)}$ over $M$
 
Moreover, the equality in (\ref{ine-i-1}) or (\ref{ine-i-2}) holds if and only if $M$ is an Einstein manifold.
\end{thm}
Observe that  Theorem \ref{thm-1} is a particular case of Theorem \ref{thm-2}   ($K=0$ ).
After the work of De Lellis and Topping,   in the case of dimension $n=3, 4$,  Y. Ge and G. Wang (\cite{GW1}, \cite{GW2}) proved that Theorem \ref{thm-1}  holds under the weaker condition of nonnegative scalar curvature. However as pointed out in  \cite{dLT},     this is surely not possible for $n\geq 5$ (it can be shown using constructions similar to the ones of Section $3$ in \cite{dLT}).  Also, Ge, Wang and Xia  \cite{GWX}  proved the case of equality in (\ref{ine-i1-1}) by a different way and gave some generalization of the De Lellis-Topping' inequalities for $k$-Einstein tensors and Lovelock curvature.

On the other hand, there is a similar phenomenon in submanifold theory.  In differential geometry, a classical  theorem states that a closed  totally umbilical surface in the Euclidean space $\mathbb{R}^3$ must be a round sphere $\mathbb{S}^2$  and its second fundamental form $A$ is a constant multiple of its metric.  This theorem is also also true for hypersurfaces in $\mathbb{R}^{n+1}$.
It is interesting to discuss the stability and rigidity of this theorem. De Lellis and M$\ddot{u}$ller \cite{dLM} obtained an optimal rigidity estimate for closed surfaces in $\mathbb{R}^3$. Recently, D. Perez \cite{P} proved the following theorem for convex hypersurfaces in $\mathbb{R}^{n+1}$:
 \begin{thm} (\cite{P})\label{thm-3}
Let $\Sigma$ be a smooth, closed and connected hypersurface in $\mathbb{R}^{n+1}, n\geq 2$ with induced Riemannian $g$ and non-negative Ricci curvature, then 
\begin{equation}\label{ine-i2}\int_{\Sigma}|A-\frac1n\overline{H}g|^2\leq \frac{n}{n-1}\int_{\Sigma}|A-\frac{H}{n}g|^2,
\end{equation}
and  equivalently 
\begin{equation}\label{ine-i-02}\int_{\Sigma}(H-\overline{H})^2\leq \frac{n}{n-1}\int_{\Sigma}|A-\frac{H}{n}g|^2,
\end{equation}
where  $A$ and  $H=\text{trace} A$ denote   the second fundamental form and the mean curvature  of $\Sigma$ respectively, $\overline{H}=\frac1{\text{Vol}_n(\Sigma)}\int_{\Sigma}H$.  In particular, the above estimate holds for smooth, closed hypersurfaces which are the boundary of a convex set in $\mathbb{R}^{n+1}$.
\end{thm}
As pointed  out by De Lellis and Topping \cite{dLT},  Perez's theorem holds even for the closed hypersufaces with nonnegative Ricci curvature when the ambient space is Einstein.   Indeed a slight modification of the  proof of Theorem  \ref{thm-3} gives

 \begin{thm}  \label{thm-4}Let $(N^{n+1}, \widetilde{g})$ be an Einstein manifold, $ n\geq 2$.
Let $\Sigma$ be a smooth, closed and connected hypersurface immersed in $N$ with  non-negative Ricci curvature, then 
\begin{equation}\label{ine-i3}\int_{\Sigma}|A-\frac1n\overline{H}g|^2\leq \frac{n}{n-1}\int_{\Sigma}|A-\frac{H}{n}g|^2,
\end{equation}
and equivalently 
\begin{equation}\label{ine-i-03}\int_{\Sigma}(H-\overline{H})^2\leq \frac{n}{n-1}\int_{\Sigma}|A-\frac{H}{n}g|^2,
\end{equation}
where $\overline{H}=\frac1{\text{Vol}_n(\Sigma)}\int_{\Sigma}H$. 

\end{thm}
Later,  Zhou and the author  (\cite{CZ}) discussed the rigidity of the equalities  in  inequalities (\ref{ine-i2}) and (\ref{ine-i-02}) and proved the following

\begin{thm}\label{thm-5} (\cite{CZ})  Let $\Sigma$ be a smooth, connected, oriented and closed  hypersurface immersed in  the Euclidean space $\mathbb{R}^{n+1}, n\geq 2,$ with  non-negative Ricci curvature. Then, 
\begin{equation}\label{eq-i-1} \int_{\Sigma}|A-\frac1n\overline{H}g|^2=\frac{n}{n-1}\int_{\Sigma}|A-\frac{H}{n}g|^2,
\end{equation}
and  equivalently
\begin{equation}\label{eq-i-01}\int_{\Sigma}(H-\overline{H})^2=\frac{n}{n-1}\int_{\Sigma}|A-\frac{H}{n}g|^2,
\end{equation}
 holds if and only if $\Sigma$ is a totally umbilical hypersurface, where $\overline{H}=\frac1{\text{Vol}_n(\Sigma)}\int_{\Sigma}H$, that is, $\Sigma$ is a distance sphere $S^n$ in $\mathbb{R}^{n+1}$.

\end{thm}
In \cite{CZ}, the authors also studied the general case for hypersurfaces without hypothesis on convexity (that is, $A\geq 0$, which is equivalent to $\text{Ric}\geq 0$ when $\Sigma$ is a closed hypersurface in $\mathbb{R}^{n+1}$).  Precisely,  the following theorem was proved.
 \begin{thm} (\cite{CZ}) \label{thm-6} Let $(N^{n+1}, \widetilde{g})$ be an Einstein manifold, $ n\geq 2$.
Let $\Sigma$ be a smooth, connected, oriented and closed   hypersurface immersed in $N$ with induced metric $g$.  Then 
\begin{equation}\label{ine-i-4}\int_{\Sigma}|A-\frac{\overline{H}}{n} g|^2\leq \frac{n}{n-1}(1+\frac{K}{\lambda_1})\int_{\Sigma}|A-\frac{H}{n}g|^2,
\end{equation}
and equivalently
\begin{equation}\label{ine-i-04}\int_{\Sigma}(H-\overline{H})^2\leq\frac{n}{n-1}(1+\frac{nK}{\lambda_1})\int_{\Sigma}|A-\frac{H}{n}g|^2,
\end{equation}
where $\lambda_1$ is the first nonzero eigenvalue of the Laplacian operator on $\Sigma$, $K\geq 0$ is a  nonnegative constant such that the Ricci curvature of $\Sigma$  satisfies $\text{Ric} \geq -K$.

Morover,  when  $N^{n+1}$ is  the Euclidean space $\mathbb{R}^{n+1}$, the hyperbolic space $\mathbb{H}^{n+1}(-1)$ or the closed hemisphere $\mathbb{S}_{+}^{n+1}(1)$, the equality in (\ref{ine-i-4}) or (\ref{ine-i-04}) holds if and only if $\Sigma$ is  a totally umbilical hypersurface,  that is, $\Sigma$ is a distance sphere $S^n$ in $N^{n+1}$.
\end{thm}

See more details in \cite{CZ}.

From \cite{dLT}, \cite{GW1}, \cite{GW2}, \cite{GWX}, \cite{C}, \cite{P} and \cite{CZ}, we observe that the inequalities mentioned above may be generalized  to symmetric $(2,0)$ tensor fields.  Applying such  unified inequalities for symmetric $(2,0)$ tensors,  we may  obtain  new inequalities besides the inequalities in the papers mentioned  above. For this purpose, in this paper, we prove that
\begin{thm}  \label{thm-r1} Let $(M, g)$ be a closed Riemannian manifold of dimension $n, n\geq 2$. Let $T$ be a symmetric   $(2,0)$-tensor field on $M$. If  the divergence $\textrm{div} T$ of $T$ and  the trace $B=\text{tr}T$  of $T$ satisfy
$\textrm{div}T=c\nabla B,$ where $c$ is a constant, then 
\begin{equation}\label{ine-r1}(nc-1)^2\int_{M}(B-\overline{B})^2\leq n(n-1)(1+\frac{nK}{\lambda_1})\int_{M}|T-\frac{B}{n}g|^2,
\end{equation}
and equivalently,
\begin{align}\label{ine-r2} (nc-1)^2&\int_M|T-\frac{\overline{B}}{n}g|^2\nonumber\\
&\leq \left[(nc-1)^2+(n-1)(1+\frac{nK}{\lambda_1})\right]\int_{M}|T-\frac{B}{n}g|^2,
\end{align}
where    $\overline{B}=\frac{\int_M Bdv}{\text{Vol}(M)}$ denotes the average of $B$ over $M$, $\lambda_1$ denotes the first nonzero eigenvalue of Laplace operator on $M$
 and $K$ is nonnegative constant such that  the Ricci curvature of $M$ satisfies $\text{Ric}\geq -(n-1)K$.
 
Furthur, assume the Ricci curvature $\text{Ric}$ of $M$ is positive. If $c\neq \frac{1}{n}$,  the following conclusions (i), (ii) and (iii) are equivalent; if  $c=\frac 1n$, the following (i) and (ii) are equivalent.

(i)  the equality in (\ref{ine-r1}) or  in (\ref{ine-r2}) holds;

(ii)  $T=\frac{{B}}{n}g$ holds on $M$;

(iii) $T=\frac{\overline{B}}{n}g$ holds on $M$.
\end{thm}
Take $K=0$ in Theorem \ref{thm-r1}. We obtain corresponding inequalities with universal constants as follows,

\begin{thm}  \label{thm-r2} Let $(M, g)$ be a closed Riemannian manifold of dimension $n, n\geq 2,$ with nonnegative Ricci curvature. Let $T$ be a symmetric   $(2,0)$-tensor field on $M$. If  the divergence $\text{div} T$ of $T$ and  the trace $B=\text{tr}T$  of $T$ satisfy
$\text{div}T=c\nabla B,$ where $c$ is a constant, then
\begin{equation}\label{ine-r3}(nc-1)^2\int_{M}(B-\overline{B})^2\leq n(n-1)\int_{M}|T-\frac{B}{n}g|^2,
\end{equation}
and equivalently, 
\begin{equation}\label{ine-r4}(nc-1)^2\int_M|T-\frac{\overline{B}}{n}g|^2\leq \left[(nc-1)^2+1\right]\int_{M}|T-\frac{B}{n}g|^2,
\end{equation}
where   $\overline{B}=\frac{\int_M Bdv}{\text{Vol}(M)}$ denotes the average of $B$ over $M$.

Furthur, assume the Ricci curvature $\text{Ric}$ of $M$ is positive. If $c\neq \frac{1}{n}$,  the following conclusions (i), (ii) and (iii) are equivalent; if  $c=\frac 1n$, the following (i) and (ii) are equivalent.

(i)  the equality in (\ref{ine-r3}) or  in (\ref{ine-r4}) holds;

(ii)  $T=\frac{{B}}{n}g$ holds on $M$;

(iii) $T=\frac{\overline{B}}{n}g$ holds on $M$.

\end{thm}

It is a known fact that if $(M^n, g), n\geq 2,$ is a connected Riemannian manifold of dimension $n$. If $T=\frac{B}{n}g$ and $\text{div} T=c\nabla B$, where $c\neq \frac 1n$ is a constant, then $B$ is constant on $M$ and thus $T$ is constant multiple of its metric $g$ (see Proposition \ref{prop-1}). Hence Theorems \ref{thm-r1} and \ref{thm-r2} discuss the stability and rigidity of this fact for closed manifolds.
Especially, take $T=\text{Ric}, A$, etc, in  Theorems \ref{thm-r1} and \ref{thm-r2}. We obtain the corresponding inequalities mentioned beforeTheorems \ref{thm-r1}.  In this paper,  we will obtain two other applications as follows.

First we deal with   $r$th mean curvatures  of a closed hypersurface in a space form.  Let $(N^{n+1}_a, \widetilde{g})$ be an $(n+1)$-dimensional space form with constant sectional curvature $a$, $n\geq 2$. Assume  $(\Sigma, g)$ is a connected oriented closed hypersurface  immersed in $(N^{n+1}_a, \widetilde{g})$  with the induced metric $g$. Associated with the second fundmental form $A$ of $\Sigma$, we have    $r$th mean curvatures $H_r$ of $\Sigma$ and the Newton transformations $P_r, 0\leq r\leq n,$ (see their definition and related notations in Section \ref{sec-rmean}).
Since  Reilly \cite{R} introduced them, there have been much work in studying high order $r$-mean curvatures (cf. for instance,  \cite{Ro}, \cite{BC}, \cite{CR}, \cite{ALM}).  It can be verified  that if the Newton transformations $P_r$ satisfy $P_r=\frac{\text{tr}P_r}{n}g$ on $\Sigma$, $\Sigma$ has constant $r$th mean curvature and thus $P_r$ is a constant multiple of its metric $g$ (see Proposition \ref{prop-1} and Section \ref{sec-rmean}). In this paper, we discuss the stability and rigidity of this fact. 

In addition,  although it is true that  a closed  immersed totally umbilical hypersurface $\Sigma$ (that is, $\Sigma$ satisfies $P_1=\frac{\text{tr}P_1}{n}g$) in  $\mathbb{R}^{n+1}$ must be a round sphere $\mathbb{S}^n$, it is unknown, to our best knowledge,  if it is true that  a closed  immersed   hypersurface $\Sigma$ satisfying $P_r=\frac{\text{tr}P_r}{n}g$ in  $\mathbb{R}^{n+1}$ must be a round sphere $\mathbb{S}^n$ for $r\geq 2$. When  $\Sigma$ is embedded, 
Ros \cite{Ro1}, \cite{Ro2}  showed  that the round spheres are the only closed embedded hypersurfaces with constant $r$th mean curvature in  $\mathbb{R}^{n+1}$, $2\leq r\leq n$ (recall Alexandrov theorem says that  the round spheres are the only closed embedded hypersurfaces  in $\mathbb{R}^{n+1}$  with constant mean curvature \cite{A}).  Hence  the round spheres are the only closed embedded hypersurfaces  in  $\mathbb{R}^{n+1}$ with $P_r=\frac{\text{tr}P_r}{n}g$, $2\leq r\leq n$.

In Section \ref{sec-rmean}, we prove the following 
\begin{thm}  \label{thm-rm1} Let $(N_a^{n+1}, \widetilde{g})$ be a space form with constant sectional curvature $a$, $ n\geq 2$.
Assume $\Sigma$ is a smooth connected oriented closed  hypersurface immersed in $N$ with induced metric $g$.  Then  for $2\leq r\leq n$,
\begin{equation}\label{ine-rm1}(n-r)^2\int_{\Sigma}(s_r-\overline{s}_r)^2\leq n(n-1)(1+\frac{nK}{\lambda_1})\int_{\Sigma}|P_r-\frac{(n-r)s_r}{n}g|^2,
\end{equation}
and equivalently, 
\begin{equation}\label{ine-rm2}\int_{\Sigma}|P_r-\frac{(n-r)\overline{s}_r}{n}g|^2\leq n\left[1+\frac{(n-1)K}{\lambda_1}\right]\int_{\Sigma}|P_r-\frac{(n-r)s_r}{n}g|^2,
\end{equation}
where  $s_r=\text{tr}P_r= \left( \begin{smallmatrix} n\\ r \end{smallmatrix} \right)H_r$,  $\overline{s}_r=\frac{\int_M s_rdv}{\text{Vol}(M)}$ denotes the average of $s_r$ over $\Sigma$, $\lambda_1$ is the first nonzero eigenvalue of the Laplacian operator on $\Sigma$, and $K\geq 0$ is a  nonnegative constant such that the Ricci curvature of $\Sigma$  satisfies $\text{Ric} \geq -K$.
 Furthur, it holds that
 \begin{itemize}
\item  [ 1)] if the Ricci curvature $\text{Ric}$ of $\Sigma$ is positive,  the following conclusions (i), (ii), and (iii)are equivalent,

(i) the equality in (\ref{ine-rm1}) or (\ref{ine-rm2}) holds;

(ii)$P_r=\frac{{(n-r)s_r}}{n}g$ holds on $\Sigma$;

(iii) $P_r=\frac{{(n-r)\overline{s}_r}}{n}g$ holds on $\Sigma$.
\item  [2)] if $\Sigma$ is embedded in the Euclidean space $\mathbb{R}^{n+1}$  and the Ricci curvature $\text{Ric}$ of $\Sigma$ is positive,   the equality in (\ref{ine-rm1}) or (\ref{ine-rm2}) holds if and only if $\Sigma$ is a round sphere $\mathbb{S}^{n+1}$ in $\mathbb{R}^{n+1}$. 

\end{itemize}

 \end{thm}
  Take $K=0$ in Theorem \ref{thm-rm1}, we prove the following
 \begin{thm}  \label{thm-rm2} Let $(N_a^{n+1}, \widetilde{g})$ be a space form with constant sectional curvature $a$, $ n\geq 2$.
Assume $\Sigma$ is a smooth connected oriented closed  hypersurface immersed in $N$ with induced metric $g$. If $\Sigma$ has nonnegative Ricci curvature,  then  for $2\leq r\leq n$,
\begin{equation}\label{ine-rm03}(n-r)^2\int_{\Sigma}(s_r-\overline{s}_r)^2\leq n(n-1)\int_{\Sigma}|P_r-\frac{(n-r)s_r}{n}g|^2,
\end{equation}
and equivalently, 
\begin{equation}\label{ine-rm4}\int_{\Sigma}|P_r-\frac{(n-r)\overline{s}_r}{n}g|^2\leq n\int_{\Sigma}|P_r-\frac{(n-r)s_r}{n}g|^2.
\end{equation}
where  $s_r=\text{tr}P_r= \left( \begin{smallmatrix} n\\ r \end{smallmatrix} \right)H_r$,  $\overline{s}_r=\frac{\int_M s_rdv}{\text{Vol}(M)}$ denotes the average of $s_r$ over $\Sigma$. Moreover, it holds that
 \begin{itemize}
\item  [ 1)] if the Ricci curvature $\text{Ric}$ of $\Sigma$ is positive,  the following conclusions (i), (ii), and (iii)are equivalent,

(i) the equality in (\ref{ine-rm03}) or (\ref{ine-rm4}) holds;

(ii)$P_r=\frac{{(n-r)s_r}}{n}g$ holds on $\Sigma$;

(iii) $P_r=\frac{{(n-r)\overline{s}_r}}{n}g$ holds on $\Sigma$.
\item  [2)] if $\Sigma$ is embedded in the Euclidean space $\mathbb{R}^{n+1}$  and the Ricci curvature $\text{Ric}$ of $\Sigma$ is positive,  the equality in (\ref{ine-rm03}) or (\ref{ine-rm4}) holds if and only if $\Sigma$ is a round sphere $\mathbb{S}^{n+1}$ in $\mathbb{R}^{n+1}$. 

\end{itemize}

\end{thm}

Second, we consider  the $k$-scalar curvatures of  locally conformally flat closed  manifolds (see their definition in Section \ref{sec-rscalar}). The $k$-scalar curvatures  have been much studied in recent years (cf, for instance, \cite{G}, \cite{V1}, \cite{V2}, etc) since they were first introduced by Viaclovsky \cite{V1}.  When $M$ is locally conformally flat,
 we obtain an almost-Schur type lemma for $k$-scalar curvatures, $k\geq 2$, as follows
 \begin{thm} \label{thm-ks1} Let $(M^n,g)$ be an $n$-dimensional closed locally conformally flat  manifold, $n \geq 3$. Then for $2\leq k\leq n$,  the $k$-scalar curvature $\sigma_k(S_g)$ and the Newton transformation $T_k$ associated with the Schouten tensor $S_g$ satisfy
\begin{align}\label{ine-ks1} (n-k)^2&\int_{M}(\sigma_k(S_g)-\overline{\sigma}_k(S_g))^2\nonumber\\
&\leq n(n-1)(1+\frac{nK}{\lambda_1})\int_{M}|T_k-\frac{(n-k)\sigma_k(S_g)}{n}g|^2,
\end{align}
and equivalently, 
\begin{align}\label{ine-ks2}&\int_{M}|T_k-\frac{(n-k)\overline{\sigma}_k(S_g)}{n}g|^2\nonumber\\
&\leq n\left[1+\frac{(n-1)K}{\lambda_1}\right]\int_M |T_k-\frac{(n-k)\sigma_k(g)}{n}g|^2,
\end{align}
where   $\overline{\sigma}_k(S_g)=\frac{\int_M \Sigma_k(S_g) dv}{\text{Vol}(M)}$ denotes the average of $\sigma_k(S_g)$ over $M$,  $\lambda_1$ is the first nonzero eigenvalue of the Laplacian operator on $\Sigma$, $K\geq 0$ is a  nonnegative constant such that the Ricci curvature of $\Sigma$  satisfies $\text{Ric} \geq -K$.

 Moreover, if the Ricci curvature $\text{Ric}$ of $M$ is positive,  the following conclusions (i), (ii), and (III) are equivalent,
 
(i) the equality in (\ref{ine-ks1}) or (\ref{ine-ks2}) holds;

(ii)  $T_k=\frac{{\sigma_k}(S_g)}{n}g$ holds on $M$;

(iii) $T_k=\frac{{\sigma_k}(\overline{S}_g)}{n}g$ holds on $M$.
\end{thm}

Take $K=0$ in Theorem \ref{thm-ks1}, we have the following result:

\begin{thm} \label{cor-ks2} Let $(M^n,g)$ be an $n$-dimensional locally conformally flat closed Riemannian manifold with nonnegative Ricci curvature, $n \geq 3$. Then for $2\leq k\leq n$, the $k$-scalar curvature $\sigma_k(S_g)$ and the Newton transformation $T_k$ associated with the Schouten tensor $S_g$ satisfy
\begin{align}(n-k)^2&\int_{M}(\sigma_k(S_g)-\overline{\sigma}_k(S_g))^2\nonumber\\
&\leq n(n-1)\int_{M}|T_k-\frac{(n-k)\sigma_k(S_g)}{n}g|^2,
\end{align}
and equivalently, 
\begin{align}\int_{M}|T_k-\frac{(n-k)\overline{\sigma}_k(S_g)}{n}g|^2
\leq n\int_M |T_k-\frac{(n-k)\sigma_k(g)}{n}g|^2.
\end{align}
 Moreover, if the Ricci curvature $\text{Ric}$ of $M$ is positive, (i), (ii), and (III) are equivalent,
 
(i)   the equality in (\ref{ine-ks1}) or (\ref{ine-ks2}) holds;

(ii)  $T_k=\frac{{\sigma_k}(S_g)}{n}g$ holds on $M$. 

(iii) $T_k=\frac{{\sigma_k}(\overline{S}_g)}{n}g$ holds on $M$.

\end{thm}

The rest of this paper is organized as follows. In Section \ref{sec-p},  we prove Theorems \ref{thm-r1} and \ref{thm-r2}. In Section \ref{sec-newton}, we recall the definitions of Newton transformation and $r$th symmetric function associated with a symmetric endomorphism of an $n$-dimensional vector space. In Section \ref{sec-rmean},  we prove Theorem \ref{thm-rm1} by applying  \ref{thm-r1}.  In Section \ref{sec-rscalar},  we  prove Theorem \ref{thm-ks1} by applying  \ref{thm-r1}.


\section{Proof of theorems on symmetric $(2,0)$-tensors}\label{sec-p}

First we give  some notations.  Assume $(M,g)$ is an $n$-dimensional closed (that is, compact and without boundary) Riemannian manifold. Let $\nabla$ denote the Levi-Civita connection on $(M,g)$ and also the induced connections on tensor bundles on $M$.    Let $T$ denote a symmetric $(2,0)$-tensor field on $M$. $\text{tr}$ denotes the trace of a tensor. $B=\text{tr} T=T_i^{{}{i}}=g^{ij}T_{ij}$ denotes the trace of $T$.  Here and thereafter  we use Einstein summation convention.  
Denote by $\overline{B}=\frac{\int_M B}{\text{Vol}(M)}$  the average of $B$ over $M$ and $\text{\r T}=T-\frac{B}{n}g$.
Denote by  $ \text{div}$   the divergence of  tensor field.  For $T$,  $\text{div}T=\text{tr}\nabla T$ is a $(1,0)$-tensor. Under the local coordinates $\{x_i\}$ on $M$,    $\text{div}T=g^{ij}(\nabla_{i}T_{jk})dx^k$, where $\nabla_{i}T_{jk}=(\nabla_{\partial_i}  T)(\partial_j,\partial_k).$   

The following proposition is a well known fact, which was mentioned in the introduction.
\begin{prop} \label{prop-1} Assume $(M^n, g), n\geq 2,$ is a connected Riemannian manifold of dimension $n$. If $T=\frac{B}{n}g$ and $\text{div} T=c\nabla B$, where $c\neq \frac 1n$ is a constant, then $B=const$ on $M$ and $T$ is constant multiple of its metric $g$.
\end{prop}
\begin{rema} Proposition \ref{prop-1} can be proved directly by noting $T=\frac{B}{n}g$ implies   the identity  $\text{div}T=\frac{\nabla B}{n}$.
\end{rema}

 The  argument of Theorem \ref{thm-r1} is similar to the one of  Theorem \ref{thm-2} (i.e. \cite{C} Thm.1.2) and in the case of $K=0$,  the one of  Theorem \ref{thm-1} (i.e. \cite{dLT} Thm.0.1) 
 
 \bigskip

{\it Proof of Theorem \ref{thm-r1}}.  Obviously, it  suffices to prove the case $c\neq\frac 1n$.  By the assumption
$\text{div}T=c\nabla B,$
\begin{align} \text{div} \text{\r T}=\text{div} T-\text{div} (\frac{B}{n}g)=\text{div} T-\frac{\nabla B}{n}=\frac{nc-1}{n}\nabla B. \label{eq-p-1}
\end{align}
Let $f$ be the unique solution of the following Poisson equation on $M$:
\begin{equation}\label{ine-p-1}
\Delta f=B-\overline{B}, \quad\quad \int_{M}f=0.
\end{equation}
By  (\ref{eq-p-1}), (\ref{ine-p-1}) and Stokes'  formula, 
\begin{align}\label{ine-p-01}
\int_{M}(B-\overline{B})^2&=\int_{M}(B-\overline{B})\Delta f=-\int_{M}\left<\nabla B,\nabla f\right>\nonumber\\
&=-\frac{n}{nc-1}\int_{M}\left<\text{div} \text{\r T},\nabla f\right>\nonumber\\
&=\frac{n}{nc-1}\int_{M}\left< \text{\r T},\nabla^2 f\right>\nonumber\\
&=\frac{n}{nc-1}\int_{M}\left< \text{\r T},\nabla^2 f-\frac1n(\Delta f)g\right>\nonumber\\
&\leq\frac{n}{|nc-1|}\left(\int_{M}| \text{\r T}|^2\right)^{\frac12}\left[\int_{M}|\nabla^2 f-\frac1n(\Delta f)g|^2\right]^{\frac12}\nonumber\\
&=\frac{n}{|nc-1|}\left(\int_{M}| \text{\r T}|^2\right)^{\frac12}\left[\int_{M}|\nabla^2 f|^2-\frac1n\int_{M}(\Delta f)^2\right]^{\frac12}
\end{align}
Recall the Bochner formula
$$\frac12\Delta|\nabla f|^2=|\nabla^2f|^2+\text{Ric}(\nabla f,\nabla f)+\left<\nabla f,\nabla(\Delta f)\right>,$$
and integrate it. By the Stokes'  formula, we have
\begin{equation}\label{eq-p-01}\int_{M}|\nabla^2f|^2=\int_{M}(\Delta f)^2-\int_{M}\text{Ric}(\nabla f,\nabla f).
\end{equation}
By (\ref{ine-p-01}) and  (\ref{eq-p-01}), 
\begin{align}\label{ine-p-2}
&\int_{M}(B-\overline{B})^2\nonumber\\
&\qquad\leq\frac{n}{|nc-1|}\left(\int_{M}| \text{\r T}|^2\right)^{\frac12}\left[\frac{n-1}{n}\int_{M}(\Delta f)^2-\int_{M}\text{Ric}(\nabla f,\nabla f)\right]^{\frac12}.
\end{align}
By (\ref{ine-p-1}),   $f\equiv 0$ if and only if $B-\overline{B}\equiv 0$ on $M$. In this case,  (\ref{ine-r1}) and (\ref{ine-r2}) obviously hold. In the following we only consider that $f$ is not identically zero.
Since the Ricci curvature has $\text{Ric}\geq -(n-1)K$  on $M$, 
\begin{equation}\label{eq-p-001}\int_{M}\text{Ric}(\nabla f, \nabla f)\geq -(n-1)K\int_{M}|\nabla f|^2.
\end{equation}
By (\ref{eq-p-001}),  (\ref{ine-p-2}) turns into
\begin{align}\label{ine-p-3}
&\int_{M}(B-\overline{B})^2\nonumber\\
&\quad\leq\frac{n}{|nc-1|}\left(\int_{M}| \text{\r T}|^2\right)^{\frac12}\left [\frac{n-1}{n}\int_{M}(\Delta f)^2+(n-1)K\int_{M}|\nabla f|^2\right ]^{\frac12}.
\end{align}
Since the first nonzero eigenvalue  $\lambda_1$ of Laplace operator on $M$ satisfies
\[\lambda_1=\displaystyle\inf\{ \frac{\int_M |\nabla \varphi|^2}{\int_M \varphi^2};  \varphi\in C^{\infty}(M) \textrm{ is not identically zero and }\int_M \varphi=0 \},
\]   
\begin{align}
\int_{M} |\nabla f|^2 &=-\int_{M}f\Delta f=-\int_{M}f(B-\overline{B})\nonumber\\
&\leq \left(\int_{M}f^2\right)^{\frac12}\left[\int_{M}(B-\overline{B})^2\right]^{\frac12}\nonumber\\
&\leq \left(\frac{\int_{M} |\nabla f|^2}{\lambda_1}\right)^{\frac12}\left[\int_M (B-\overline{B})^2\right]^{\frac12}.\nonumber
\end{align}
Then \begin{equation}\label{ine-p-4}
\int_{M} |\nabla f|^2 \leq  \frac{1}{\lambda_1}\int_M (B-\overline{B})^2.
\end{equation}
Substitute (\ref{ine-p-4}) into (\ref{ine-p-3}) and note that  $K\geq 0$. We have
\begin{align}\label{ine-p-04}
&\int_{M}(B-\overline{B})^2\nonumber\\
&\leq\frac{n}{|nc-1|}\left(\int_{M}| \text{\r T}|^2\right)^{\frac12}\left[\frac{n-1}{n}\int_{M}(B-\overline{B})^2+\left(\frac{(n-1)K}{\lambda_1}\right)\int_M (B-\overline{B})^2\right]^{\frac12}\nonumber\\
&=\frac{n^{\frac12}(n-1)^{\frac12}}{|nc-1|}\left(1+\frac{nK}{\lambda_1}\right)^{\frac12}\left[\int_{M}| \text{\r T}|^2\right]^{\frac12}\left[\int_M (B-\overline{B})^2\right]^{\frac12}
\end{align}
(\ref{ine-p-04}) implies that
\begin{align}\int_{M}(B-\overline{B})^2\leq \frac{n(n-1)}{(nc-1)^2}(1+\frac{nK}{\lambda_1})\int_{M}| \text{\r T}|^2.
\end{align}
Thus we have inequality  (\ref{ine-r1}):
\begin{equation}(nc-1)^2\int_{M}(B-\overline{B})^2\leq n(n-1)(1+\frac{nK}{\lambda_1})\int_{M}|T-\frac{B}{n}g|^2.\nonumber
\end{equation}
By the  identity
$|\text{T}-\frac{\overline{B}}{n}g|^2=|T-\frac{B}{n}g|^2+\frac{1}{n}(B-\overline{B})^2,$
we have inequality (\ref{ine-r2}):
\begin{equation*} (nc-1)^2\int_M|T-\frac{\overline{B}}{n}g|^2\leq \left[(nc-1)^2+(n-1)(1+\frac{nK}{\lambda_1})\right]\int_{M}|T-\frac{B}{n}g|^2.
\end{equation*}

Now with the assumption of positivity of  Ricci curvature $\text{Ric}$ of $M$, we may prove the case of equalities in  (\ref{ine-r1}) and (\ref{ine-r2}). Obviously,  if $T=\frac{{B}}{n}g$ holds on $M$,  the equalities in (\ref{ine-r1}) and (\ref{ine-r2}) hold.
On the other hand,  suppose  the equality in (\ref{ine-r1}) (or equivalently (\ref{ine-r2})) holds. If $c=\frac 1n$, it is obvious that $T=\frac{{B}}{n}g$ on $M$.  If $c\neq \frac 1n$, we may  take $K=0$. By the proof of (\ref{ine-r1}),  the equality  in (\ref{ine-r1}) holds if and only if 
\begin{itemize}
\item [1)] $\text{Ric}(\nabla f, \nabla f)=0$ on $M$ and 
\item [2)] $T-\frac{B}{n}g$ and $\nabla^2 f-\frac1n(\Delta f)g$ are linearly dependent.
\end{itemize}
Note  $\text{Ric}>0$ and 1). It must holds that $\nabla f\equiv 0$ on $M$.  Then $f\equiv 0$. Thus $B=\overline{B}$ on $M$. By  (\ref{ine-r1}), we obtain that $T=\frac{{B}}{n}g$ on $M$. Hence conclusions (i) and (ii) are equivalent. Obviously (iii) implies (ii). When $c\neq \frac1n$, if (ii) holds, by the above argument, (ii) implies $B=\overline{B}$ on $M$. Thus (iii) also holds.

\qed 

We have a corollary of Theorem \ref{thm-r1} as follows:
 \begin{cor}\label {cor-1} Let $(M^n, g), n\geq 2,$ be a closed Riemannian manifold of dimension $n$. Let $T$ be a symmetric   $(2,0)$-tensor field on $M$. If  the divergence $\text{div} T$ of $T$ and  the trace $B=\text{tr}T$  of $T$ satisfy
$\text{div}T=c\nabla B,$ where $c\neq \frac1n$ is a constant, then  
\begin{equation}\label{ine-c1}\int_{M}(B-\overline{B})^2\leq C_{(Kd^2)}\int_{M}|T-\frac{B}{n}g|^2,
\end{equation}
and 
\begin{equation}\label{ine-c2}\int_{M}|T-\frac{\overline{B}}{n}g|^2\leq \overline{C}_{(Kd^2)}\int_{M}|T-\frac{B}{n}g|^2,
\end{equation}
where  $K$ is a positive constant. such that  the  Ricci curvature of $M$ satisfies $\text{Ric}\geq -(n-1)K,$   $d$ denotes the diameter of $M$ and $C_{(Kd^2)}$ and $ \overline{C}_{(Kd^2)}$ are constants only depending on $Kd^2$.
\end{cor}

{\it Proof of Corollary \ref{cor-1}.}  When $\text{Ric}\geq -(n-1)K$,  where constant$K>0$,
Li and Yau \cite{LY} proved   that  the first nonzero eigenvalue $\lambda_1$ has the lower bound: 
$$\lambda_1\geq \alpha=\frac{1}{(n-1)d^2\exp [1+\sqrt{1+4(n-1)^2Kd^2}]},$$ 
where $d$ denotes the diameter of $M$.
So $$\frac{K}{\lambda_1}\geq \frac{K}{\alpha}=(n-1)Kd^2\exp [1+\sqrt{1+4(n-1)^2Kd^2}].$$
By Theorem \ref{thm-r1}, we obtain  inequality  (\ref{ine-c1}) with the constant 
$$C_{(Kd^2)})=\frac{4n(n-1)}{(n-2)^2}\left(1+n(n-1)Kd^2\exp [1+\sqrt{1+4(n-1)^2Kd^2}]\right).$$
Inequality  (\ref{ine-c1}) implies inequality  (\ref{ine-c2}).
\qed

\begin{rema}There are other  lower estimates $\alpha$ of $\lambda_1$ using the diameter $d$ and negative lower bound $-(n-1)K$ of Ricci curvature (for example, see \cite{KMYZ}). Hence we may have other values of constant $C_{(Kd^2)}$ and $ \overline{C}_{(Kd^2)}$.
\end{rema}

\section{Newton transformations and the $r$th elementary symmetric function}\label{sec-newton}
 Let $\sigma_r:\mathbb{R}^r\rightarrow \mathbb{R}$ denote the elementary symmetric function in $\mathbb{R}^n$ given by
$$\sigma_r(x_{i_1},\cdots,x_{i_r})=\displaystyle\sum_{i_1<\cdots<i_r}x_{i_1}\ldots x_{i_r}, 1\leq r\leq n.$$
Let $V$ be an $n$-dimensional vector space and  $A:  V\rightarrow V$ be a symmetric linear transformation.  If $\eta_1, \cdots, \eta_n$ are
the eigenvalues of $A$  corresponding the orthonormal eigenvectors $\{e_i\}, i=1,\ldots, n$ respectively, 
define the $r$th symmetric functions $\sigma_r(A)$  associated with $A$ by 
\begin{eqnarray}
\sigma_0(A)&=&1,\nonumber
\\ \sigma_r(A)&=&\sigma_r(\eta_{i_1},\ldots, \eta_{i_k}), 1\leq r\leq n.
\end{eqnarray}
For convenience of the notation, we simply denote $\sigma_r(A)$ by $\sigma_r$ if there is  no  confusion.
The Newton
transformations $P_r: V\rightarrow V, $ associated with $A, 0\leq r\leq n,$  are defined by
$$P_0=I,$$
$$P_r=\displaystyle\sum_{j=0}^r(-1)^j\sigma_{r-j}A^j=\sigma_rI-\sigma_{r-1}A+...+(-1)^rA^r, r=1,\ldots,n.$$
By definition, 
$P_r=\sigma_rI-AP_{r-1},  P_n=0.$ 
 It was proved in  \cite{R} that  $P_r$ has the following basic properties:
\begin{eqnarray}
&(i)&P_r(e_i)=\frac{\partial \sigma_{r+1}}{\partial \eta_i}e_i;\nonumber
\\&(ii)&\text{tr}(P_r)
=(n-r)\sigma_r;\nonumber
\\&(iii)&\text{tr}(AP_r)
=(r+1)\sigma_{r+1}.\nonumber
\end{eqnarray}
Obviously, each $P_r$ corresponds a symmetric $(2,0)$-tensor on $V$, still denoted by $P_r$.

\section{High order  mean curvatures of hypersurfaces in space forms}\label{sec-rmean}

Assume $(N, \widetilde{g})$ is an $(n+1)$-dimensional Riemannian manifold, $n\geq 2$. Suppose 
 $(\Sigma, g)$ is a smooth connected oriented closed hypersurface immersed in  $(N, \widetilde{g})$ with induced metric $g$.   Let  $\nu$ denote the outward unit normal to $\Sigma$ and $A=(h_{ij})$ denote the second fundamental form $A: T_p\Sigma\otimes_s T_p\Sigma\rightarrow \mathbb{R}$, defined by   $A(X,Y)=-\left<\widetilde{\nabla}_XY,\nu\right>$, where  $X, Y\in T_p\Sigma, p\in \Sigma$, $ \widetilde{\nabla}$ denote the Levi-Civita connection of $(N,\widetilde{g})$.  $A$ determines an equivalent $(1,1)$-tensor, called   the shape operator $A$ of $\Sigma$:  $T_p\Sigma \rightarrow T_p\Sigma$, given by $AX=\widetilde{\nabla}_X\nu$.   $\Sigma$ is called totally umbilical if $A$ is multiple of its metric $g$ at every point of $\Sigma$, that is, $A=\frac{\text{tr}A}{n}g$ on $\Sigma$. Now we recall the definition of $r$th mean curvatures of a hypersurface, which was introduced by Reilly \cite{R}, (cf.\cite{Ro}).

Let  $\eta_i, i=1,\ldots, n$ denote the principle curvatures  of $\Sigma$ at $p$, which are the eigenvalues of $A$ at $p$ corresponding the orthonormal eigenvectors $\{e_i\}, i=1.\ldots, n$ respectively.  By Section \ref{sec-newton}, we have 
 the $r$th symmetric functions $\sigma_r(A)$ associated with $A$,  denoted by $s_r=\sigma_r(A)$ and   the Newton
transformations $P_r$ associated with $A$ at $p$, $0\leq r\leq n$.  We have 
\begin{Def} The $r$th mean
curvature $H_r$ of $\Sigma$ at $p$ is defined by
$s_r=\left( \begin{matrix} n\\ r \end{matrix}\right) H_r,$ $ 0\leq r\leq n.$
\end{Def}
For instance, $H_1=\frac{s_1}{n}=\frac{H}{n}$  (in this paper, we also call $H=\text{tr}A$ the mean curvature of $\Sigma$, conformal to   the  previous related papers  \cite{P}, \cite{CZ}, etc).  
$H_n$ is the
Gauss-Kronecker curvature. 
When the ambient space $N$ is a space form $N_a^{n+1}$ with constant sectional curvature $a$,
\begin{align*}&\text{Ric}=(n-1)aI+HA-A^2.\\
&R=\text{tr}\text{Ric}=n(n-1)c+H^2-|A|^2=n(n-1)a+2s_2.
\end{align*}
Hence $H_2$ is, modulo a constant, the
scalar curvature of $\Sigma$. 

One of  the known properties of $P_r$   is the following 
 \begin{lem} \label{lem-r-1}( \cite{R}, cf \cite{Ro}, or  \cite{ALM}) When the ambient space is a space form $N_a^{n+1}$, $\text{div} P_r=0, 0\leq r\leq n.$
\end{lem}
Now we prove Theorem \ref{thm-rm1}.

{\it Proof of Theorem \ref{thm-rm1}.} By Section \ref{sec-newton}, $\text{tr}P_r=(n-r)s_r.$ 
 Denote by $\overline{s}_r=\frac{\int_{\Sigma}s_r}{\text{Vol}(\Sigma)}$. By Lemma \ref{lem-r-1},  $\text{div} P_r=0$. Take $T=P_r$ and $B=(n-r)s_r$ in Theorem \ref{thm-r1}.  We have 
\begin{equation}\label{ine-rm1-1}(n-r)^2\int_{\Sigma}(s_r-\overline{s}_r)^2\leq n(n-1)(1+\frac{nK}{\lambda_1})\int_{\Sigma}|P_r-\frac{(n-r)s_r}{n}g|^2,\nonumber
\end{equation}
and equivalently, 
\begin{equation}\label{ine-rm2-2}\int_{\Sigma}|P_r-\frac{(n-r)\overline{s}_r}{n}g|^2\leq n\left(1+\frac{(n-1)K}{\lambda_1}\right)\int_{\Sigma}|P_r-\frac{(n-r)s_r}{n}g|^2,\nonumber
\end{equation}
which are (\ref{ine-rm1}) and (\ref{ine-rm2}) respectively.

Now we prove conclusions 1) and 2) in Theorem \ref{thm-rm1}.  If the Ricci curvature of $\Sigma$ is positive, by Theorem \ref{thm-r1},  conclusion 1)  holds and $s_r=\overline{s}_r$ is constant on $\Sigma$.
If  $\Sigma$ is also embedded in $\mathbb{R}^{n+1}$,  by the Ros' theorem \cite{Ro2} that a closed embedded hypersurface in $\mathbb{R}^{n+1}$ with constant $r$th mean curvature must be a distance sphere $\mathbb{S}^{n+1}$, $2\leq r\leq n$, we obtain conclusion 2).

\qed

\begin{rema}  If $r=1$, $P_1=s_1I-A=HI-A.$ $P_1$ is equivalent to  the symmetric $(2,0)$-tensor $P_1=Hg-A$. So (\ref{ine-rm1})  turns to
\begin{align}
\int_{\Sigma}\label{ine-rm5} (H-\overline{H})^2&\leq \frac{n}{n-1}(1+\frac{nK}{\lambda_1})\int_{\Sigma}|Hg-A-\frac{(n-1)H}{n}g|^2\nonumber \\
&= \frac{n}{n-1}(1+\frac{nK}{\lambda_1})\int_{\Sigma}|A-\frac{H}{n}g|^2.
\end{align}
In particular, if  $K=0$, 
\begin{align}
\int_{\Sigma}\label{ine-rm6} (H-\overline{H})^2 
\leq \frac{n}{n-1}\int_{\Sigma}|A-\frac{H}{n}g|^2.
\end{align}
(\ref{ine-rm5}) and  (\ref{ine-rm6}) are  (\ref{ine-i-04}) and (\ref{ine-i-03}) respectively, which were proved in \cite{CZ} and   \cite{P} respectively  if $\Sigma$ is a closed hypersurface immersed in an Einstein manifold. This is because $\text{div}P_1=0$ even if the ambient space is Einstein.

When $r=2$, we have $2s_2=R-n(n-1)a$, 

$P_2=s_2I-s_1A+s_0A^2=\frac{R-(n-2)(n-1)a}{2}I-\text{Ric}$, and by direct computation, 
\begin{align}
P_2-\frac{(n-2)s_2}{n}g=\frac{R}{n}I-\text{Ric}.\nonumber
\end{align}
As a symmetric $(2,0)$-tensor, $P_2=\frac{R}{n}g-\text{Ric}.$
Hence (\ref{ine-rm1}) turns to
\begin{equation}\int_{\Sigma}(s_2-\overline{s}_2)^2\leq \frac{n(n-1)}{(n-2)^2}(1+\frac{nK}{\lambda_1})\int_{\Sigma}|P_2-\frac{(n-2)s_2}{n}g|^2,\nonumber
\end{equation}
which is 
\begin{equation}\label{ine-rm3}\int_{\Sigma}(R-\overline{R})^2\leq \frac{4n(n-1)}{(n-2)^2}(1+\frac{nK}{\lambda_1})\int_{\Sigma}|\text{Ric}-\frac{R}{n}g|^2.
\end{equation}
(\ref{ine-rm3}) was proved in \cite{C} and in the case of $K=0$, was proved in \cite{dLT}.

If $r=n$, (\ref{ine-rm1})  is trivial.
\end{rema}

\section{$k$-scalar curvature of locally conformal flat manifolds. }\label{sec-rscalar}

We first  recall the definition of the k-scalar curvatures of a Riemannian manifold,  introduced by Viaclovsky in \cite{V1}. If $(M^n,g)$ be an $n$-dimensional Riemannian manifold, $n \geq 3$, the
Schouten tensor of $M$ is $$S_g=\frac{1}{n-2}\left(Ric-\frac{1}{2(n-1)}Rg\right).$$
By definition, $S_g: TM\rightarrow TM$ is a symmetric $(1,1)$-tensor  field. By Section \ref{sec-newton}, we have the symmetric $k$th function  $\sigma_k(S_g)$ and the Newton transformations $T_k(S_g)=T_k$ associated with $S_g$, $1\leq k\leq n$.  We call $\sigma_k(S_g)$  the k-scalar curvatures of $M$ 

It was proved that
\begin{lem}\label{lem-k1} (\cite{V1}) If $(M,g)$ is locally conformally flat, then for $1\leq k\leq n$,  $div T_k(S_g)=0.$
\end{lem}
Because of Lemma \ref{lem-k1},  we can applying Theorem \ref{thm-r1} to $T_k(S_g)$ to  obtain Theorem \ref{thm-ks1}. 
\begin{rema} When $k=1$,  $\sigma_1(S_g)=\text{tr}S_g=\frac{R}{2(n-1)}$,  $T_1=\sigma_1(S_g)I-S_g$. As a symmetric $(2,0)$-tensor, $T_1=-\frac{1}{n-2}(\text{Ric}-\frac{Rg}{2}).$ Hence (\ref{ine-ks1}) turns to (\ref{ine-i-1})
\begin{equation}\int_{M}(R-\overline{R})^2\leq \frac{4n(n-1)}{(n-2)^2}\left(1+\frac{nK}{\lambda_1}\right)\int_{M}|\text{Ric}-\frac{R}{n}g|^2,\nonumber
\end{equation}
and in particular, if $K=0$, (\ref{ine-ks1})  turns to (\ref{ine-i1-1})
\begin{equation}\int_{M}(R-\overline{R})^2\leq \frac{4n(n-1)}{(n-2)^2}\int_{M}|\text{Ric}-\frac{R}{n}g|^2.\nonumber
\end{equation}
 (\ref{ine-i-1}) and (\ref{ine-i1-1}) were proved in \cite{C} and \cite{dLT} respectively without the hypothesis that $M$ is locally conformally flat.  The reason is that $\text{div}T_1=0$ (the contracted second Bianchi identity) holds on any Riemannian manifold.
\end{rema}

\bigskip
\bigskip
\noindent  Xu Cheng\\Insitituto de Matem\'atica\\Universidade
Federal Fluminense - UFF\\Centro, Niter\'{o}i, RJ 24020-140 Brazil
\\e-mail:xcheng@impa.br

\bigskip


\begin{thebibliography}{9999}


\bibitem[A] {A}  
A.D. Alexandrov, Uniqueness theorems for surfaces in the large V, Vestnik Leningrad Univ. Math. 13 (1958), 5Ð8; English translation: AMS Transl. 21 (1962), 412Ð416.

\bibitem[ALM] {ALM}
L. Al\'ias, J. Lira  and J. Malacarne,  Constant higher-order mean curvature hyper- surfaces in Riemannian spaces,  Jour. Inst. Math. Jussieu, 5 (2006), 527Ð562. 

\bibitem[BC] {BC}  J.L. Barbosa  and A.G. Colares,  Stability of hypersurfaces with constant r-mean curvature. Ann Global Anal and Geom 15: 277Ð297,1997.

\bibitem[C]{C} X.  Cheng,   A generalization of almost-Schur lemma for closed Riemannian manifolds, Ann  Global Anal and Geom, (28 June 2012), 1-8, (online).

\bibitem[CR]{CR} X. Cheng  and H. Rosenberg,  Embedded positive constant r-mean curvature hypersurfaces in $M^m \times \mathbb{R}$, Anais da Academia Brasileira de Ci\^encias (Annals of the Brazilian Academy of Sciences) (2005) 77(2), 183 -199.


\bibitem[CZ]{CZ} X.  Cheng and D. Zhou, Rigidity for nearly umbilical hypersurfaces in space forms,  arXiv:1208.1786, 2012.

\bibitem[dLM]{dLM} C. De Lellis and S. M$\ddot{u}$ller, Optimal rigidity estimates for nearly umbilical surfaces. J. Differential Geom. 69 (2005) 75-110.

\bibitem[dLT]{dLT} C. De Lellis and P. Topping, Almost -Schur Lemma,  Calc. Var. and PDE, 43 (2012) 347--354;  arXiv:1003.3527v2 [math.DG] 7 May 2011.


\bibitem[GW1]{GW1} Y. Ge and G. Wang,   An almost Schur theorem on $4$-dimensional manifolds, Proc. Amer. Math. Soc. 140 (2012), 1041-1044.

\bibitem[GW2]{GW2} Y. Ge and G. Wang,   A new conformal invariant on $3$-dimensional manifolds, arXiv:1103.3838, 2011.

\bibitem[GWX]{GWX} Y. Ge,  G. Wang and Chao Xia,   On problems related to an inequality of De Lellis and Topping (preprint), 2011. 

\bibitem[G]{G}  P. Guan, Topics	 in Geometric Fully Nonlinear Equations,	Lecture Notes, http://www.math.mcgill.ca/guan/notes.html

\bibitem[dLM]{dLM} C. De Lellis and S. M$\ddot{u}$ller, Optimal rigidity estimates for nearly umbilical surfaces. J. Differential Geom. 69 (2005) 75-110.

\bibitem[dLT]{dLT} C. De Lellis and P. Topping, Almost -Schur Lemma,  Calc. Var. and PDE, 43 (2012) 347--354;  arXiv:1003.3527v2 [math.DG] 7 May 2011.

\bibitem[LY]{LY} P. Li and S.T. Yau, Eigenvalues of a compact Riemannian manifold, AMS Proc. Symp. Pure Math., 36 (1980), 205-239.

\bibitem[KMYZ]{KMYZ} M. Kalka, E. Mann,  D.Yang and A. Zinger, The Exponential Decay Rate of the Lower Bound for the First Eigenvalue of Compact Manifolds, International Journal of Mathematics (IJM), Volume: 8, Issue: 3(1997) pp. 345-355.                                                                                     


\bibitem[P]{P} D. Perez, On nearly umbilical hypersurfaces, thesis, 2011.

\bibitem[R]{R} R. C. Reilly, Variational properties of functions of the mean curvatures for hypersur- faces in space forms, J. Differential Geometry, 8 (1973), 465-477.

\bibitem[Ro1]{Ro1}A.  Ros, Compact hypersurfaces with constant scalar curvature and a congruence theorem. J. Differential Geom. 27 (1988), 215-220.
\bibitem[Ro2]{Ro2} A. Ros,  Compact hypersurfaces with constant higher order mean curva tures. Rev. Mat. Iberoamericana 3 (1987), 447-453.

\bibitem[Ro]{Ro} H. Rosenberg, Hypersurfaces of constant curvature in space forms, Bull. Sci. Math., 117 (1993), 211-239.

\bibitem[V1]{V1}  J. Viaclovsky, Conformal geometry, contact geometry, and the calculus of variations, Duke Math. J.,
101 (2000), 283-316.

\bibitem[V2]{V2} J. Viaclovsky, Conformal geometry and fully nonlinear equations, Inspired by S. S. Chern, 435Ð460,
Nankai Tracts Math. 11 World Sci. Publ., Hackensack, NJ, 2006
\end{thebibliography}
\end{document}